# Riemann's Zeta Function.
# $\xi(s)$ via analytic Continuation
# inside an infinite Loop.

Renaat Van Malderen

December 2013

*Abstract*

The paper deals with the analytic entire function $\xi(s)$ closely related to Riemann's Zeta Function $\zeta(s)$. A formula is obtained for $\xi(s)$ essentially within the so-called critical strip. This is achieved by applying Cauchy's integral formula to an infinite loop encircling the critical strip. In the obtained formula a remarkable role is played by a special type of complex function known as "Incomplete Gamma Function". Numerical examples verifying the obtained formula are included.

## 1. Introduction.

The present paper deals with entire functions f(z) i.e. functions analytic over the entire complex plane. Consider a simple closed contour C in this plane and a function f(z) as defined above. Suppose we have some type of analytic expression for f(z) outside and on the contour C, but we lack a "convenient" corresponding expression for f(z) in part or the whole inside of C.

Cauchy's integral formula along the contour C allows to determine f(z) for any point $z_0$ inside C:

$$f(z_0) = \frac{1}{2\pi i} \oint_c \frac{f(z)dz}{z - z_0} \quad (1)$$

(1) in fact represents the analytic continuation of $f(z)$ from the contour C to its interior.

In case we are able to solve (1), it might result in some "convenient" expression for $f(z_0)$ within C.

An obvious example (although not an entire function) would be Riemann's Zeta function $\zeta(s)$ (We are here using the customary complex variable s= σ+it. As will become clear further on, we will also work with the variable z= x + it where s= z+1/2). Simple expressions for $\zeta(s)$ exist in the range σ>1, and via its "functional equation" also for the range σ<0. In the so-called critical strip i.e. 0≤σ≤1, the situation is different. Here a number of expressions exist for $\zeta(s)$ based on either integrals, series expansions or combinations there-of, but these are for analytic purposes in general not very insightful, or for computational purposes maybe inconvenient. It might therefore be of interest to look for further alternative expressions. Although $\zeta(s)$ itself is not directly suitable for the approach described in this paper, its well-known related entire function [1,p.16]:

$$\xi(s) = \frac{s(s-1)}{2} \Gamma\left(\frac{s}{2}\right) \pi^{-\frac{s}{2}} \zeta(s) \quad (2)$$

suits the situation well as explained further on.

## 2. Further restrictions on f(z).

In addition to f(z)=f(x+it) assumed to be entire, we also require it to be even and real on the x-axis. Taking these restrictions into account, the power series expansion of f(z) around the origin equals:

$$f(z) = \sum_{n=0}^{\infty} a_{2n} z^{2n} \quad (3)$$

with all $a_{2n}$ real.

(3) meets the requirements spelled out above:
f(-z) = f(z) and f(x)=real.

(3) also implies that f(z) is real on the t-axis. Indeed putting z=ρe$^{iφ}$:

$$f(z) = \sum_{n=0}^{\infty} a_{2n} z^{2ni\varphi} \rho^{2n} \text{ and } f(\bar{z}) = \sum_{n=0}^{\infty} a_{2n} z^{-2ni\varphi} \rho^{2n} \quad (4)$$

(4) implies $f(\bar{z}) = \overline{f(z)} = f(-\bar{z})$ (5).
These latter symmetry properties include f(it) being real.

## 3. Transition from a finite to an infinite contour.

Consider a contour C in the z plane, oriented in counter clockwise direction and consisting of the four segments:

C1: z = a+it with a>0 and –T ≤ t ≤ T with T > 0.
C2: z = x+iT with –a ≤ x ≤ a.
C3: z = -a+it with –T ≤ t ≤ T.
C4: z = x-iT with –a ≤ x ≤ a.

Given an entire function f(z), f(z$_0$) for a point z$_0$ inside C is given by (1). The integral around C may be split-up into four parts I$_1$, I$_2$, I$_3$, I$_4$ corresponding to segments C$_1$, C$_2$, C$_3$, C$_4$.

$$\oint_C = I_1 + I_2 + I_3 + I_4 \quad (6)$$

We now let C change in size by keeping its width constant but increasing its vertical size 2T indefinitely. Expression (1) will remain valid during this stretching. In addition to the restrictions put on f(z) in section-2, we make two more assumptions:

A1: For the improper integrals $I_1$ and $I_3$:
$\lim_{\pm T \to \pm\infty} I_1$ and $\lim_{\pm T \to \pm\infty} I_3$ converge each independently in both the positive and negative t direction.

A2: With a being kept constant:
$\lim_{T \to \infty} I_2 = 0$ and $\lim_{-T \to -\infty} I_4 = 0$

Under these conditions (6) reduces to $I_1 + I_3$ and becomes:

$$f(z_0) = \frac{1}{2\pi i}\int_{-\infty}^{\infty} \frac{f(a+it)i\,dt}{a+it-z_0} - \frac{1}{2\pi i}\int_{-\infty}^{\infty} \frac{f(-a+it)i\,dt}{-a+it-z_0} \quad (7)$$

From (5) we know f(-a+it)=f(a-it):

$$f(z_0) = \frac{1}{2\pi}\int_{-\infty}^{\infty} \frac{f(a+it)\,dt}{a+it-z_0} + \frac{1}{2\pi}\int_{-\infty}^{\infty} \frac{f(a-it)\,dt}{a-it+z_0} \quad (8)$$

In the particular case $z_0 = it_0$, i.e. on the t-axis:

$$f(it_0) = \frac{1}{2\pi}\int_{-\infty}^{\infty} dt \left[\frac{f(a+it)}{a+i(t-t_0)} + \frac{f(a-it)}{a-i(t-t_0)}\right] \quad (9)$$

### 4. An archetypical example: Cosh z.

We determine the value of Cosh z on the vertical z=it using (9):

$$\text{Ch}(it_0) = \frac{1}{2\pi}\int_{-\infty}^{\infty} dt \left[\frac{\text{Ch}(a+it)}{a+i(t-t_0)} + \frac{\text{Ch}(a-it)}{a-i(t-t_0)}\right]$$

$$\text{Ch}(it_0) = \frac{1}{2\pi}\int_{-\infty}^{\infty} dt \left[\frac{e^a e^{it} + e^{-a}e^{-it}}{2(a+i(t-t_0))} + \frac{e^a e^{-it} + e^{-a}e^{+it}}{2(a-i(t-t_0))}\right]$$

Substituting $\tau = t - t_0$:

$$\text{Ch}(it_0) = \frac{1}{2\pi}\int_{-\infty}^{\infty} d\tau \left[\frac{e^{a+it_0}e^{i\tau} + e^{-a-it_0}e^{-i\tau}}{2(a+i\tau))} + \frac{e^{a-it_0}e^{-i\tau} + e^{-a+it_0}e^{i\tau}}{2(a-i\tau)}\right] \quad (10)$$

The integrals in (10) are convergent. Using well-known formulas [2,p269] from complex analysis (residue calculus) with a>0:

$$\frac{1}{2\pi}\int_{-\infty}^{\infty} \frac{e^{i\tau}d\tau}{a+i\tau} = e^{-a}$$

$$\frac{1}{2\pi} \int_{-\infty}^{\infty} \frac{e^{-i\tau} d\tau}{a - i\tau} = e^{-a} \qquad (11)$$

$$\frac{1}{2\pi} \int_{-\infty}^{\infty} \frac{e^{-i\tau} d\tau}{a + i\tau} = 0$$

$$\frac{1}{2\pi} \int_{-\infty}^{\infty} \frac{e^{i\tau} d\tau}{a - i\tau} = 0$$

(10) yields:

$\text{Ch}(it_0) = \frac{e^{it_0} + e^{-it_0}}{2} = \cos t_0$ as was to be expected.

### 5. Behaviour of ξ(s) for large |t| Values.

$$\xi(s) = \frac{s(s-1)}{2} \pi^{-\frac{s}{2}} \Gamma\left(\frac{s}{2}\right) \zeta(s) \quad (12)$$

As explained in [1, p.17] ξ(s) is an even function with respect to s=1/2, i.e. $\xi\left(s - \frac{1}{2}\right) = \xi\left(\frac{1}{2} - s\right)$ or in terms of z=s-1/2: $\xi(z) = \xi(-z)$. Also ξ(z) is real on both the x and t axis. To avoid confusion between s and z, where deemed helpful, we will indicate explicitly which variable is being used. It is our intent to obtain $\xi(z)$ on the critical strip by integrating formula (8) along the line s=2+it. For s= 2+it we have:

$$\left|\Gamma\left(\frac{2+it}{2}\right)\right| \sim \sqrt{2\pi} \left(\frac{t}{2}\right)^{1/2} e^{\left(\frac{-\pi|t|}{4}\right)} \text{ see}[3,\text{p.60}]$$

Also: $|\zeta(2 + it)| \leq \zeta(2) = \frac{\pi^2}{6}$
From (12) eventually:

$$|\xi(s = 2 + it)| < K|t|^{5/2} e^{\left(\frac{-\pi|t|}{2}\right)} \quad (13) \text{ with } K \text{ some positive constant.}$$

### 6. The case of ξ(s): Steps towards solving the integral of formula (8).

As already indicated above, in our further calculations we will pick a particular value for the parameter a in (7), (8), (9): a=3/2. Notice this is in terms of z; the line s=2+it becomes $z = \frac{3}{2} + it$ in z terms.
We work out the detailed expression for ξ(s=2+it):

$$\xi(s) = \frac{s}{2}(s - 1)\pi^{-\frac{s}{2}} \Gamma\left(\frac{s}{2}\right) \zeta(s) = (s - 1)\Gamma\left(1 + \frac{s}{2}\right) \pi^{-\frac{s}{2}} \zeta(s) \quad (14)$$

We want a linear t-term in the denominator to avoid problems of convergence of integrals in our manipulations further down:

$$\xi(s) = \frac{s-1}{1+\frac{s}{2}} \Gamma\left(2+\frac{s}{2}\right) \pi^{-\frac{s}{2}} \zeta(s)$$

With $s=2+it$

$$\xi(s = 2 + it) = \frac{2(1+it)}{\pi(4+it)} \Gamma\left(3+\frac{it}{2}\right) \pi^{-\frac{it}{2}} \zeta(2+it) \quad \text{with}$$

$$\Gamma\left(3+\frac{it}{2}\right) = \int_0^\infty e^{-\lambda} \lambda^{2+\frac{it}{2}} d\lambda \qquad (15)$$

$$\zeta(2+it) = \sum_{n=1}^\infty \frac{e^{-\left(\frac{it\ln(n^2)}{2}\right)}}{n^2} \qquad (16)$$

$$\xi(s = 2 + it) = \frac{2}{\pi}\frac{(1+it)}{(4+it)}\pi^{-it/2} \sum_{n=1}^\infty \frac{e^{-\left(\frac{it\ln(n^2)}{2}\right)}}{n^2} \int_0^\infty e^{-\lambda}\lambda^{2+it/2} d\lambda \qquad (17)$$

Consider the first of the integrals of (8) and insert $\xi(s = 2 + it)$ (which is of course the same as $\xi(z = 3/2 + it)$) from (17) into it:

$$\frac{1}{2\pi}\frac{2}{\pi} \int_{-\infty}^{+\infty} dt \frac{\xi(z = \frac{3}{2} + it)}{\frac{3}{2} + it - z_0}$$

$$= \frac{1}{\pi^2} \int_{-\infty}^{+\infty} dt \frac{(1+it)\pi^{-it/2}}{(4+it)(\frac{3}{2}+it-z_0)} \sum_{n=1}^\infty \frac{e^{-\left(\frac{it\ln(n^2)}{2}\right)}}{n^2} \int_0^\infty e^{-\lambda}\lambda^{2+it/2} d\lambda \qquad (18)$$

Since $\zeta(s) = \sum_{n=1}^\infty n^{-s}$ converges absolutely and uniformly in the half plane σ>1+ε with ε>0 [4, p. 251], we can move the summation sign in (18) to the front:

$$(18) = \frac{1}{\pi^2} \sum_{n=1}^\infty \frac{1}{n^2} \int_{-\infty}^{+\infty} dt \int_0^\infty d\lambda \frac{(1+it)e^{-\left(\frac{it\ln(n^2)}{2}\right)}e^{-\lambda}\lambda^{2+it/2}}{(4+it)(\frac{3}{2}+it-z_0)} \qquad (19)$$

The function in two variables t and λ in the range -∞<t<+∞, 0<λ<+∞ which figures in the double integral (19) is continuous in both variables combined (For proper definition see e.g. [5, p. 119]). Moreover the integral representing the Gamma function $\int_0^\infty e^{-\lambda}\lambda^{2+it/2} d\lambda$ converges uniformly and absolutely for the range considered [3, p. 51]. As such we are allowed to change the order of integration [5, p. 314] and after combining exponents:

$$(19) = \frac{1}{\pi^2} \sum_{n=1}^\infty \frac{1}{n^2} \int_0^\infty d\lambda e^{-\lambda}\lambda^2 \int_{-\infty}^{+\infty} dt \frac{(1+it)e^{\left(\frac{it}{2}\ln\left(\frac{\lambda}{\pi n^2}\right)\right)}}{(4+it)(\frac{3}{2}+it-z_0)} \qquad (20)$$

The arguments leading to (20) are equally valid for the second integral in (8).

### 7. $\xi(s)$ on the Critical Line $s= \frac{1}{2} +it_0$

From (9) and using the conclusions of the previous section in terms of $z_0=it_0$:

$$\xi(it_0) = \frac{1}{\pi^2}\sum_{n=1}^{\infty}\frac{1}{n^2}\int_0^{\infty} d\lambda e^{-\lambda}\lambda^2(I_1 + I_2) \quad (21)$$

with:

$$I_1 = \int_{-\infty}^{+\infty} dt \, \frac{(1+it)e^{\left(\frac{it\beta}{2}\right)}}{(4+it)\left(\frac{3}{2}+i(t-t_0)\right)}$$

and

$$I_2 = \int_{-\infty}^{+\infty} dt \, \frac{(1-it)e^{\left(-\frac{it\beta}{2}\right)}}{(4-it)\left(\frac{3}{2}-i(t-t_0)\right)}$$

$\quad (22)$

where $\beta = \ln\left(\frac{\lambda}{\pi n^2}\right)$

Partial fraction expansion yields:

$$I_1 = \int_{-\infty}^{+\infty} dt \, e^{\left(\frac{it\beta}{2}\right)}\left[\frac{A}{(4+it)} + \frac{B}{\left(\frac{3}{2}-it_0+it\right)}\right]$$

$$I_2 = \int_{-\infty}^{+\infty} dt \, e^{\left(-\frac{it\beta}{2}\right)}\left[\frac{C}{(4-it)} + \frac{D}{\left(\frac{3}{2}+it_0-it\right)}\right]$$

with:

$$A = \frac{3}{\frac{5}{2}+it_0} \qquad B = \frac{-\frac{1}{2}+it_0}{\frac{5}{2}+it_0}$$

$$C = \frac{3}{\frac{5}{2}-it_0} \qquad D = \frac{-\frac{1}{2}-it_0}{\frac{5}{2}-it_0}$$

We use residue calculus and Jordan's theorem [2, p.219] to solve $I_1$ and $I_2$. t in this context is treated as a complex variable. Table -1 shows the location of the involved poles and corresponding residues.

For $I_1$ its two poles occur in the upper semi-plane. To apply Jordan's theorem we need to consider the contour consisting of $-\infty < t < +\infty$ and the upper semi-circle at infinity on which the contribution to $I_1$ will be zero for $\beta > 0$ which requires $\ln\left(\frac{\lambda}{\pi n^2}\right) > 0$ or $\lambda > \pi n^2$. For $\lambda < \pi n^2$, $I_1$ equals zero. For $I_2$ the situation is the opposite: The sign in the exponent of $I_2$ is negative but since its poles are in the lower semi-plane, we still need $\beta > 0$, so the same condition applies as for $I_1$. The rule for adding residues assumes an anti-clockwise direction of the contour, but the direction of $I_2$ runs from left to right. So we need to change signs of both $R_3$ and $R_4$. Taking all this into account we obtain [2, p.207]:

$$I_1 + I_2 = 2\pi i \sum \text{Residues }(R_1, R_2, -R_3, -R_4)\, u(\lambda - \pi n^2) \quad (23)$$

where $u(\lambda - \pi n^2)$ stands for the unit step function which equals zero for $\lambda < \pi n^2$ and 1 for $\lambda > \pi n^2$. This changes the lower limit of integration in (21) from zero to $\pi n^2$.

Table -1

|  | Term | Pole location | Residue |
|---|---|---|---|
| $I_1$ | $\dfrac{1}{4+it}$ | $t_1 = 4i$ | $R_1 = \dfrac{-3iE(-2\beta)}{\frac{5}{2}+it_0}$ |
|  | $\dfrac{1}{\frac{3}{2}-it_0+it}$ | $t_2 = t_0 + \frac{3}{2}i$ | $R_2 = \dfrac{i\left(\frac{1}{2}-it_0\right)E\left(-\frac{\beta}{2}\left(\frac{3}{2}-it_0\right)\right)}{\frac{5}{2}+it_0}$ |
| $I_2$ | $\dfrac{1}{4-it}$ | $t_3 = -4i$ | $R_3 = \dfrac{3iE(-2\beta)}{\frac{5}{2}-it_0}$ |
|  | $\dfrac{1}{\frac{3}{2}+it_0-it}$ | $t_4 = t_0 - \frac{3}{2}i$ | $R_4 = \dfrac{-i\left(\frac{1}{2}+it_0\right)E\left(-\frac{\beta}{2}\left(\frac{3}{2}+it_0\right)\right)}{\frac{5}{2}-it_0}$ |

Adding up the residues in (23) with appropriate signs:

$$I_1 + I_2 = 2\pi\left[3T_1 E(-2\beta) - T_2 E\left(-\frac{\beta}{2}\left(\frac{3}{2}-it_0\right)\right) - \overline{T_2}E\left(-\frac{\beta}{2}\left(\frac{3}{2}+it_0\right)\right)\right]u(\lambda - \pi n^2) \quad (24)$$

with:

$$\left.\begin{array}{l}
T_1 = \dfrac{5}{\frac{25}{4}+t_0^2} \quad \text{and} \quad T_2 = \dfrac{\frac{1}{2}-it_0}{\frac{5}{2}+it_0} \\[1em]
E(-2\beta) = \dfrac{\pi^2 n^4}{\lambda^2} \quad \text{and} \quad E\left(\dfrac{-3\beta}{4}\right) = \dfrac{\pi^{3/4}n^{3/2}}{\lambda^{3/4}} \\[1em]
E\left(\dfrac{it_0\beta}{2}\right) = \left(\dfrac{\lambda}{\pi n^2}\right)^{\frac{it_0}{2}} \quad \text{and} \quad E\left(-\dfrac{it_0\beta}{2}\right) = \left(\dfrac{\lambda}{\pi n^2}\right)^{-\frac{it_0}{2}}
\end{array}\right\} \quad (26)$$

Plugging (24) into (21), using (26) and after the smoke clears:

$$\xi(it_0) = 2 \sum_{n=1}^{\infty} \left\{ 3T_1 \pi n^2 \int_{\pi n^2}^{\infty} d\lambda e^{-\lambda} \right.$$
$$\left. - \frac{T_2(\pi n^2)^{\frac{-it_0}{2}}}{\pi^{\frac{1}{4}}\sqrt{n}} \int_{\pi n^2}^{\infty} d\lambda e^{-\lambda} \lambda^{\frac{5}{4}+\frac{it_0}{2}} - \frac{\overline{T_2}(\pi n^2)^{\frac{+it_0}{2}}}{\pi^{\frac{1}{4}}\sqrt{n}} \int_{\pi n^2}^{\infty} d\lambda e^{-\lambda} \lambda^{\frac{5}{4}-\frac{it_0}{2}} \right\} \quad (27)$$

An integral of the type $\int_0^\alpha e^{-\lambda}\lambda^z d\lambda = \gamma(z+1,\alpha)$ is known under the name of "Lower Incomplete Gamma function" [6, p.260], [9]. Its complement, i.e. $\int_\alpha^\infty e^{-\lambda}\lambda^z d\lambda$ is denoted by $\Gamma(z+1,\alpha)$ and is called the "Upper Incomplete Gamma function". Obviously $\gamma(z+1,\alpha) + \Gamma(z+1,\alpha) = \Gamma(z+1)$. We recognize the second and third integral in (27) as upper incomplete Gamma functions. Accordingly:

$$\xi(it_0) = 2 \sum_{n=1}^{\infty} \left\{ \frac{15\pi n^2 e^{-\pi n^2}}{\frac{25}{4}+t_0^2} - \frac{(1-2it_0)(\pi n^2)^{\frac{-it_0}{2}}}{(5+2it_0)\pi^{\frac{1}{4}}\sqrt{n}} \Gamma\left(\frac{9}{4}+\frac{it_0}{2},\pi n^2\right) \right.$$
$$\left. - \frac{(1+2it_0)(\pi n^2)^{\frac{+it_0}{2}}}{(5-2it_0)\pi^{\frac{1}{4}}\sqrt{n}} \Gamma\left(\frac{9}{4}-\frac{it_0}{2},\pi n^2\right) \right\} \quad (28)$$

Since the second and the third term in (28) are complex conjugates, the total expression (28) is real as it should be.

### 8. Alternate Form for $\xi(it_0)$.

By combining complex conjugates, (28) may be written as:

$$\xi(it_0) = 2 \sum_{n=1}^{\infty} \left\{ \frac{60\pi n^2 e^{-\pi n^2}}{25+4t_0^2} - \frac{2}{\pi^{\frac{1}{4}}\sqrt{n}} \int_{\pi n^2}^{\infty} e^{-\lambda}\lambda^{\frac{5}{4}} d\lambda \left[ A\cos\left(\frac{t_0\beta}{2}\right) + B\sin\left(\frac{t_0\beta}{2}\right) \right] \right\} \quad (29)$$

with $A = (5-4t_0^2)/(25+4t_0^2)$ and $B = 12t_0/(25+4t_0^2)$

$$\beta = \ln\left(\frac{\lambda}{\pi n^2}\right)$$

For numerical evaluation (28) is more useful than (29) since with the former, series expansions for the incomplete gamma functions may be used (see sec.10).

### 9. $\xi(s=\sigma+it)$ on the Strip $0\leq\sigma\leq1$

We will now generalize (28) for points $z_0=x_0+it_0$ off the critical line, i.e. for $|x_0|\leq 1/2$. Starting from (8) with $a=3/2$:

$$\xi(z_0) = \frac{1}{2\pi}\left[\int_{-\infty}^{\infty}\frac{dt\,\xi(s=2+it)}{\frac{3}{2}-x_0+i(t-t_0)} + \int_{-\infty}^{\infty}\frac{dt\,\xi(s=2-it)}{\frac{3}{2}+x_0-i(t-t_0)} + \right]$$

The integrals corresponding to $I_1$ and $I_2$ (see(22)) are:

$$I_1^* = \int_{-\infty}^{+\infty}dt\,\frac{(1+it)e^{\left(\frac{it\beta}{2}\right)}}{(4+it)\left(\frac{3}{2}-x_0+i(t-t_0)\right)}$$

$$I_2^* = \int_{-\infty}^{+\infty}dt\,\frac{(1-it)e^{\left(\frac{-it\beta}{2}\right)}}{(4-it)\left(\frac{3}{2}+x_0-i(t-t_0)\right)}$$

Proceeding as before, using partial fractions and identifying poles and residues we have:

$$\frac{(1+it)}{(4+it)\left(\frac{3}{2}-x_0+i(t-t_0)\right)} = \frac{A}{(4+it)} + \frac{B}{\left(\frac{3}{2}-x_0+i(t-t_0)\right)}$$

with $A = \frac{3}{\frac{5}{2}+x_0+it_0}$ and $B = \frac{-\frac{1}{2}+x_0+it_0}{\frac{5}{2}+x_0+it_0}$

$$\frac{(1-it)}{(4-it)\left(\frac{3}{2}+x_0-i(t-t_0)\right)} = \frac{C}{(4-it)} + \frac{D}{\left(\frac{3}{2}+x_0-i(t-t_0)\right)}$$

with $C = \frac{3}{\frac{5}{2}-x_0-it_0}$ and $D = \frac{-\left(\frac{1}{2}+x_0+it_0\right)}{\frac{5}{2}-x_0-it_0}$

Table-2 below shows the resulting poles and residues.

| Poles | Residues |
|---|---|
| $I_1^*$ Term $1/(4+it)$ | |
| $t_1^* = 4i$ | $R_1^* = \dfrac{-3iE(-2\beta)}{\frac{5}{2}+x_0+it_0}$ |
| $I_1^*$ Term $1/(\frac{3}{2}-x_0+i(t-t_0))$ | |
| $t_2^* = t_0 + i(\frac{3}{2}-x_0)$ | $R_2^* = \dfrac{i(\frac{1}{2}-x_0-it_0)E(-\frac{\beta}{2}(\frac{3}{2}-x_0-it_0))}{\frac{5}{2}+x_0-it_0}$ |
| $I_2^*$ Term $1/(4-it)$ | |
| $t_3^* = -4i$ | $R_3^* = \dfrac{3iE(-2\beta)}{\dfrac{5}{2} - x_0 - it_0}$ |
| $I_2^*$ Term $1/(\frac{3}{2} + x_0 - i(t - t_0))$ | |
| $t_4^* = t_0 - i(\dfrac{3}{2} + x_0)$ | $R_4^* = \dfrac{-i(\frac{1}{2} + x_0 + it_0)E(-\frac{\beta}{2}(\frac{3}{2} + x_0 + it_0))}{\dfrac{5}{2} - x_0 - it_0}$ |

$I_1^* + I_2^* = 2\pi i \sum \text{Residues } (R_1^*, R_2^*, -R_3^*, -R_4^*)$

Continuing along similar lines as before we find for "the strip" area $0 \leq |x_0| \leq 1/2$:

$$\xi(x_0 + it_0) = 2 \sum_{n=1}^{\infty} \left\{ \frac{15\pi n^2 e^{-\pi n^2}}{\frac{25}{4} + t_0^2 - x_0^2 - 2it_0 x_0} - \frac{(\pi n^2)^{\frac{-(x_0+it_0)}{2}} \left(\frac{1}{2} - x_0 - it_0\right)}{\pi^{\frac{1}{4}}\sqrt{n}\left(\frac{5}{2} + x_0 + it_0\right)} \Gamma(\frac{9}{4} + \frac{x_0}{2} \right.$$

$$\left. + \frac{it_0}{2}, \pi n^2) - \frac{(\pi n^2)^{\frac{(x_0+it_0)}{2}} \left(\frac{1}{2} + x_0 + it_0\right)}{\pi^{\frac{1}{4}}\sqrt{n}\left(\frac{5}{2} - x_0 - it_0\right)} \Gamma(\frac{9}{4} - \frac{x_0}{2} - \frac{it_0}{2}, \pi n^2) \right\} \quad (30)$$

Obviously, in contrast to (28), formula (30) results generally in complex values.

With $z_0 = x_0 + it_0$, (30) allows to easily verify that $\xi(z_0) = \xi(-z_0)$. Using the abbreviation $\Psi(y) = \sum_{n=1}^{\infty} e^{-\pi n^2 y}$ we have $\Psi'(y) = -\sum_{n=1}^{\infty} \pi n^2 e^{-\pi n^2 y}$, $\Psi''(y) = \sum_{n=1}^{\infty} \pi^2 n^4 e^{-\pi n^2 y}$ and $\sum_{n=1}^{\infty} \pi n^2 e^{-\pi n^2} = \int_1^{\infty} \Psi''(y) dy = -\Psi'(1)$.

Based on these expressions, after some manipulation and partial integration, (30) may be transformed into:

$$\xi(z_0) = -4\Psi'(1) + (\frac{1}{2} - z_0) \int_1^{\infty} \Psi'(y) y^{\frac{1}{4} + \frac{z_0}{2}} dy + (\frac{1}{2} + z_0) \int_1^{\infty} \Psi'(y) y^{\frac{1}{4} - \frac{z_0}{2}} dy \quad (30a)$$

Some further algebra and a single step of partial integration of Riemann's formula (2) [1, p.17] shows this latter to be equivalent to (30a). Using $\frac{1}{2} + \Psi(1) + 4\Psi'(1) = 0$ (see [1,p.17]), (30a) may be shown to be equal to the formula given on [1, p.16] for $\xi(s)$.

## 10. Incomplete Gamma Functions.

Per definition :

$$\gamma(z + 1, \alpha) = \int_0^{\alpha} e^{-\lambda} \lambda^z d\lambda \quad (lower) \quad (31)$$

$$\Gamma(z + 1, \alpha) = \int_{\alpha}^{\infty} e^{-\lambda} \lambda^z d\lambda \quad (upper) \quad (32)$$

Obviously :

$$\gamma(z + 1, \alpha) + \Gamma(z + 1, \alpha) = \Gamma(z + 1) \quad (33)$$

Let $z = \beta + ik$ with $\beta$ and $k$ real and for our specific purposes we assume:

$$\frac{5}{4} - \frac{1}{4} \leq \beta \leq \frac{5}{4} + \frac{1}{4} \quad \text{or} \quad 1 \leq \beta \leq \frac{3}{2} \qquad (34)$$

which is the range of $\beta$ occurring in (28) and (30).

In these formulas $\Gamma(z + 1, \alpha)$ occurs in an infinite series as $\Gamma(\beta + 1 + ik, \pi n^2)$ with $k=t_0/2$, and one would expect that with increasing $n=1,2,3,...$, $|\Gamma(z + 1, \alpha)|$ will somehow start to decrease.
We may quantify this by putting a "crude" bound on it:

$$|\Gamma(\beta + 1 + ik, \alpha)| = \left|\int_\alpha^\infty e^{-\lambda}\lambda^{\beta+ik} d\lambda\right| \leq \left|\int_\alpha^\infty e^{-\lambda}\lambda^{\beta} d\lambda\right| \qquad (35)$$

In the range $\alpha = \pi n^2 \leq \lambda$ the ratio $\frac{e^{-\lambda}\lambda^2}{e^{-\lambda}\lambda^\beta} = \lambda^{2-\beta} > 1$.

So $\int_\alpha^\infty e^{-\lambda}\lambda^\beta d\lambda < \int_\alpha^\infty e^{-\lambda}\lambda^2 d\lambda < 2\alpha^2 e^{-\alpha}$ with $\alpha = \pi n^2$, and our "crude" bound becomes:
$$|\Gamma(\beta + 1 + ik, \pi n^2)| < 2\pi^2 n^4 e^{-\pi n^2} \qquad (36)$$

We call this bound "crude" because with $k \neq 0$ and $|k|$ increasing, the value of $|\Gamma(\beta + 1 + ik, \pi n^2|$ will decrease accordingly.

Table -3 gives the value of (36) for $n=1,2,3,4$.

Table -3

| n | $2\pi^2 n^4 e^{-\pi n^2}$ |
|---|---|
| 1 | 0.8530 |
| 2 | 0.0011 |
| 3 | 8.4E(-10) |
| 4 | 7.5E(-19) |

*Evaluating $\Gamma(z + 1, \alpha)$ numerically.*

For $z = \beta + ik$ with $\beta > 0$, expanding $\Gamma(z + 1, \alpha)$ in some kind of series does not look promising. Rather we will expand $\gamma(z + 1, \alpha)$ in a series and obtain $\Gamma(z + 1, \alpha)$ as $\Gamma(z + 1) - \gamma(z + 1, \alpha)$.
Repeated partial integration of (31) yields,

$$\gamma(z + 1, \alpha) = e^{-\alpha}\alpha^z \sum_{j=1}^{m} \frac{\alpha^j}{\prod_{r=1}^{j}(z+r)} + R_m$$

Where

$$R_m = \int_0^\alpha \frac{e^{-\lambda}\lambda^{z+m} d\lambda}{\prod_{r=1}^{j}(z+r)} = \int_0^\alpha e^{-\lambda}\lambda^z \varphi(\lambda, z, m) d\lambda \text{ with } \varphi(\lambda, z, m) = \prod_{r=1}^{m} \frac{\lambda}{(z+r)}$$

$\varphi(\lambda, z, m)$ obtains its maximum value for $\lambda=\alpha$.

$$|\varphi(\lambda, z, m)| < |\varphi(\alpha, z, m)| = \prod_{r=1}^{m} \frac{\alpha}{|z+r|} \qquad (37)$$

For any $\varepsilon>0$ we can select a minimum $m$ such that (37)<$\varepsilon$. In other words $lim_{m\to\infty} R_m = 0$. The terms of the remaining series $\gamma(z+1,\alpha) = e^{-\alpha}\alpha^z \sum_{j=1}^{\infty} \frac{\alpha^j}{\prod_{r=1}^{j}(z+r)}$ (38) represent a nul-sequence [7,p.17] and the series (38) is (absolutely) convergent as is easily established. (38) may also be written as:

$$\gamma(z+1,\alpha) = e^{-\alpha}\alpha^z \Gamma(z+1) \sum_{j=1}^{\infty} \frac{\alpha^j}{\Gamma(z+1+j)} \quad (39)$$

The terms of (38) reach their maximum (absolute) value for:

$\frac{\pi^j n^{2j}}{\prod_{r=1}^{j}|(z+r)|} \cong \frac{\pi^{j+1} n^{2j+2}}{\prod_{r=1}^{j+1}|(z+r)|}$ or $1 \cong \frac{\pi\ n^2}{|(z+j+1)|}$ or $[k^2 + (\beta+1+j)^2]^{1/2} \cong \pi n^2$. So the maximum size term is reached at about: $j_{max} \cong (\pi^2 n^4 - k^2)^{1/2} - (\beta+1)$.

Some observations:

- For $k \geq \pi n^2$, the first term is the largest.
- Starting with $k=0$ and then increasing its value, the maximum size term shifts to lower $j$ values until it reaches $j=1$.
- Let us consider the ratio of two successive terms:
$$Q_m = \frac{\pi n^2}{|z+1+j|}$$

  For $j_{max}>1$, starting with $j=1$ and then increasing $j$, $Q_j$ will initially go up until $j \cong j_{max}$. For $j>j_{max}$, $Q_j$ will monotonically decrease.
- For a given $k$, to obtain a certain accuracy in evaluating (38) i.e. $\gamma(z+1, \pi n^2)$ the $j_{max}$ value will increase with increasing $n$, and we will have to take more terms into account.

Regarding accuracy, the evaluation of (28) and (30) needs to consider several aspects:
  a. The accuracy of the results of (38) may be checked by increasing the number of terms calculated and observing whether there is no significant change anymore by an increase in terms.
  b. Comparison of the magnitude of $\gamma(z+1, \pi n^2)$ for n=1,2,3,… allows to decide at which point one can ignore further n values. The bound (36) is here also helpful.
  c. For increasing $k$ values, the value $\Gamma(z+1)$ goes quickly down to very small values and for increasing $n$ we are confronted with subtracting increasingly small complex numbers, i.e. $\Gamma(z+1) - \gamma(z+1, \pi n^2)$. We should in this respect not ignore that Stirling's formula is only an approximation, even if a very good one.
  d. In all of this we also need to take into account the accuracy achievable with available computational resources.

## 11. Numerical Results

Below, two numerical verifications of formula (28) for $\xi(it_0)$ are given and compared with the computation via formula (2). The required values for Gamma and Zeta functions were based on Stirling, [8], and formula (38).

Case 1: $t_0=0$ (in s- terms: $s=1/2$)

This is the simplest possible case.

*Using formula (2)*:

The elements required are:

$\Gamma(s/2)= \Gamma(1/4)= 3.6256099082$ (see [6, p. 255])
$\zeta(s=1/2)= -1.4603545$ (see[8])
Formula (2) then yields: $\xi(s=1/2)= 0.49712077$

*Using formula (28)*:

$\frac{24}{5}\sum_{n=1}^{5} \pi n^2 e^{-\pi e^2} = 0.65186088$   ($n=5$ provides plenty of accuracy)

$\Gamma(9/4)= 1.13300306$
$\Gamma(9/4, \pi)= 0.2574490078$
$\Gamma(9/4, 4\pi)= 0.000090845$

For the second and third term of (28) only two values $n=1,2$ were necessary to be considered. These Gamma function values were obtained via formula (38) and the difference $\Gamma(z+1) - \gamma(z+1, \alpha)$. Also $\Gamma(9/4)$ was obtained as $\Gamma(9/4, 9\pi)$ as it turned out to be more accurate than Stirling for small real values of the argument. The number of terms used in (38) was 60. The total for the second and third term in (28) equaled 0.15474008. The end-result via (28) equaled $\xi(1/2)=0.49712080$. So the difference between both results is about 3E(-8).

Case 2: $t_0=12$ (in s- terms: $s=1/2+12i$)

*Using formula (2)*:

Using Stirling $\Gamma(s/2)=\Gamma(0.25 +6i)=-0.000044668-0.000121314i$.
Using [8] $\zeta(s)=\zeta(0.5+12i)=1.015940-0.745105i$.
The resulting $\xi(0.5+12i)= 0.008823638+2.8E(-8)i$.
The last result not being precisely real is not surprising. The nature of furmula (2) means that only infinite precision in both Stirling and calculations could yield an exact real result.

*Using formula (28)*:

$\Gamma(9/4+6i)= 0.00268590+0.00392738i$ (Stirling)
$\Gamma(9/4+6i, \pi)=-0.04138484596+0.098818783i$
$\Gamma(9/4+6i, 4\pi)=-0.00008249+0.000005083i$

The last two results were obtained via (38) using respectively 85 and 100 terms and comparing the results. The term $\Gamma(9/4+6i, 9\pi)$ was already down in the range E(-11) and was dropped. The end result of (28) equaled:
$\xi(s=0.5+12i)=0.008823639$.

(28) of course provides a strictly real result and differs only in the ninth decimal with (2). No attempt was made to put rigorous bounds on all of the above. To do this would be near impossible since the involved calculations were often bordering the accuracy limits of the computing tool used (Ti-89 Titanium). This is also the reason why a moderate $t_0=12$ was chosen. With higher values one quickly ends up facing very small complex numbers.

Renaat Van Malderen
Address: Maxlaan 21, B-2640 Mortsel, Belgium
The author can be reached by email: hans.van.malderen@telenet.be